\documentclass[12pt]{article}
\usepackage{amssymb,latexsym, amsmath, amscd, array, graphicx}
\usepackage{amssymb}

\usepackage{setspace}

\usepackage{circledsteps}

\input xy
\xyoption{all}
\makeatletter 

\usepackage{setspace}
\usepackage{enumitem}

\date{}

\newtheorem{theorem}{Theorem}[section]

\newtheorem{proposition}[theorem]{Proposition}

\newcommand{\z}{{\Bbb Z}}
\newcommand{\q}{{\Bbb Q}}
\newcommand{\re}{{\Bbb R}}

\newcommand{\lo}{\rightarrow}

\newcommand{\black}{{\blacksquare}}

\newcommand{\diam}{{\rm diam}}

\newcommand{\mdim}{{\rm mdim}}
\newcommand{\mesh}{{\rm mesh}}
\newcommand{\ord}{{\rm ord}}
 \newcommand{\dd}{\textcolor{red}{\dagger}}

\begin{document}

\title{\bf Finite-to-one equivariant  maps  and \\ mean dimension}

\author{  Michael  Levin
\footnote{The author
 was supported by the Israel Science Foundation grant No. 2196/20 and 
 the Max Planck Institute for  Mathematics (Bonn) where 
 this note  was written  during his visit in 2023.}}

\maketitle

\begin{abstract}
We show that  a minimal dynamical system $(X, \z)$  on a compact metric $X$ with
 $\mdim X =d$    admits for every natural $k>d$ an equivariant map to the shift $( [0,1]^k)^\z$
  such that each fiber  of this map contains
 at most $[\frac{k}{k-d}]\frac{k}{k-d}$ points.  
\\\\
{\bf Keywords:}    Mean Dimension, Topological Dynamics
\\
{\bf Math. Subj. Class.:}  37B05 (54F45)
\end{abstract}
\begin{section}{Introduction}
A classical result in Dimension Theory \cite[p. 124]{kuratowski} says that if a compact metric $X$
is  with $\dim X = d$ then  for almost every map $f : X \lo [0,1]^k$, $k >d$,
the fibers  of $f$ contain at most $\frac{k}{k-d}$ points. In an attempt to generalize 
this result  to Mean Dimension we prove:

\begin{theorem}
\label{theorem}
 Let $(X, \z)$ be a minimal dynamical system on a compact metric $X$ with
 $\mdim X =d$  and let $k>d$ be a natural number. Then   almost every map $f : X \lo [0,1]^k$
 induces the map $f^\z : X \lo  ([0,1]^k)^\z$ whose fibers  contain
 at most $[\frac{k}{k-d}]\frac{k}{k-d}$ points. 

\end{theorem}
Here  $f^\z$ denotes  the equivariant  map
to the shift  $([0,1]^k)^\z$ defined by  $f^\z (x) =(f(x+z))_{z \in \z}$
 for $ x \in X$. In this note for an additive group $G$ acting on  $X$ we use  the notation $x+g$
 rather than  $gx$  to denote the action of   $g \in G $  on $x \in X$.
  This seems to be  more consistent  with the fact  that $0$ is the identity element of $G$.
 
Although  the right estimate for the size of the fibers of $f^\z$ in Theorem \ref{theorem}
should  probably be  $[\frac{k}{k-d}]$, Theorem \ref{theorem} still implies that
 $f^\z$ is an embedding  for $k >2d$ (the case of a minimal action 
of the Gutman-Tsukamoto theorem \cite{gutman-tsukamoto})  and  $f^\z$
is  finite-to-one  for $k=[d]+1$ (the case that motivated this note).

The proof of Theorem \ref{theorem} is self-contained and, basically, an interplay
of two beautiful  constructions in Topology, namely, Kolmogorov's covers (used in Kolmogorov's
proof of the continuous version of Hilbert's 13th problem \cite{kolmogorov}) 
and Lindenstrauss'
 level functions (of Rokhlin-type 
 used in Lindenstrauss' proof of the embedding theorem
\cite{lindenstrauss}). 

Another essential ingredient of the proof of Theorem \ref{theorem} 
(that to a certain extent bridges between  Kolmogorov's covers and Lindenstrauss' level functions)  is  
Borel's construction, aka  the mapping torus  in Topological Dynamics, and the following
property.

\begin{theorem}
\label{borel}
For any dynamical system $(X,\z)$ on a compact metric $X$  one has  \\
$\mdim X \times_\z \re=\mdim X$ where $X \times_\z \re$  is  Borel's construction for $(X,\z)$.
\end{theorem}
The case of a minimal action on $X$  in Theorem \ref{borel} was obtained by 
Gutman and Jin \cite{gutman} based on
 Lindenstrauss'  characterization of mean dimension via 
 metric mean dimension \cite{lindenstrauss}.  
  We prove Theorem \ref{borel}
 by a direct  argument which does not require any restriction  on the action. Theorem \ref{borel} together
 with Kolmogorov's covers and Lindenstrauss' level functions seem to  provide (at least
 for minimal actions) an elementary  replacement of the 
 signal processing technique used in the proof of Gutman-Tsukamoto's
 embedding theorem \cite{gutman-tsukamoto}.
 
 The note is organized as follows: 
 Section \ref{notations} - a description of the notations used in the note,
 Section \ref{review-mean-dimension} - a review of mean dimension,
 Section \ref{review-borel} - a review of Borel's  construction,
 Section \ref{proof-borel} -    a proof of Theorem \ref{borel}, 
 Section  \ref{review-kolmogorov} - a review of Kolmogorov's covers,
 Section \ref{review-lindenstrauss} - a review of Lindenstrauss' level functions,
 Section \ref{proof-theorem} - a proof of Theorem \ref{theorem}.

\end{section}
\begin{section}{Notations} 
\label{notations}
Everywhere in this note $i,j,m,n,l$ and $k$ stand for natural numbers (=non-negative integers)
and $z$ for an integer number (although in most cases   it is stated explicitly).

For an additive $G$ group acting on $X$ we use 
the notation $x +g$ rather than  $gx$ to denote  the action of $g \in G$ on $x \in X$ and
the expression $x+(-g)$ shortens  to $x-g$.
  We  also denote
$A +g=\{ x +g : x \in A \}$ for $A\subset X$, ${\cal A}+g =\{A +g : A \in {\cal A}\}$
for a collection $\cal A$ of subsets of $X$,  $ A+H=\{ x+g: x\in A,  g \in H \}$
for  $H\subset G$ and ${\cal A}+H=\{A +H : A \in {\cal A}\}$.

For a function $f : X \lo Y$ and collections ${\cal A}$ and ${\cal B}$ of subsets
of $X$ and $Y$ respectively we denote
$f({\cal A})=\{ f(A): A \in {\cal A} \}$ and  $f^{-1}({\cal B})=\{f^{-1}(B): B \in {\cal B}\}$.

A map means a continuous function.
For a finite open cover $\cal U$ of a compact metric $X$ the nerve of $\cal U$ is denoted by $N({\cal U})$
(a finite simplicial complex with $\ord {\cal U}=\dim N({\cal U})+1$) and a canonical map
$f : X \lo N({\cal U})$  means a map constructed on the basis of  a partition of unity
subordinate to ${\cal U}$  and an important  property of $f$  is that the fibers of $f$  refine ${\cal U}$.

\end{section}
\begin{section}{ Mean  Dimension}
\label{review-mean-dimension}
Let $(X,\z)$ be a dynamical system on a compact metric $X$ and $d$ a positive real number. 
Let us write  $\mdim X \leq d^-$ if for every $\epsilon >0$ there is a finite closed
cover $\cal A$  of $X$ such that $\mesh ({\cal A}+z) < \epsilon$ for every
 integer $z\in \z$
satisfying $0\leq zd < \ord {\cal A}$.  We will write $\mdim X > d^-$ if
$\mdim  X \leq d^-$ does not hold. Note that slightly enlarging the elements of
$\cal A$ to open sets  we may assume that in the above definition
 that $\cal A$ is an open cover. We leave to the reader to verify
that Gromov's mean dimension $\mdim X$ can be defined as
 $\mdim X =\inf \{ d : \mdim X \leq d^- \}$ (and
 $\mdim X =\infty$ if the defining set is empty). Note that $1^- < \mdim [0,1]^\z =1$.

 Let ${\cal A}_0, \dots, {\cal A}_n$ be   covers of $X$. 
 The notation ${\cal A}_0 \vee \dots \vee {\cal A}_n$ is
 used to denote the cover 
 $\{ A_0 \cap \dots \cap A_n : A_i \in {\cal A}_i, 0\leq i \leq n\}$. 
 Assume that each ${\cal A}_i$ is an open finite cover of $X$ and
 recall that for a canonical map  $f_i : X \lo K_i=N({\cal A}_i)$  to the simplicial complex
 $K_i$ which is the nerve of ${\cal A}_i$ one has that the fibers of 
 $f_i$ refine the elements of ${\cal A}_i$ and $\dim K_i  +1 ={\ord {\cal A}_i} $.
 Thus for $f=(f_0, \dots, f_n) : X \lo  K_0 \times \dots \times K_n$
  we have that the fibers of $f$ refine ${\cal A}=
  {\cal A}_0 \vee \dots \vee {\cal A}_n$
  and therefore  ${\cal A}$ can be refined by an open cover of  $X$ of 
  $\ord \leq \dim K  +1=\dim K_0 + \dots + \dim K_n +1$. 
  
  Also note  that if for  (not necessarily open) covers ${\cal A}_i$ of $X$
  we have that  $\mesh ({\cal A}_i +z)< \epsilon $ for every $0\leq i \leq n$ and
  $0\leq z < m, z \in \z$,  then for the cover
  ${\cal B}={\cal A}_0 \vee ({\cal A}_1-z_1 )\vee \dots \vee ( {\cal A}_n -z_n)$
  where $z_i=im$
  we have that $\mesh ( {\cal B}+z) < \epsilon $ for every $0\leq z < m(n+1) , z \in \z$.

\end{section}

\begin{section}{Borel's  construction}
\label{review-borel} Let $\z$ act on a compact metric $X$.
Consider the induced action of $\z$ on $X \times \re$ defined by
$(x,t)+z=(x+z, t+z), z\in\z$.  Borel's construction $X \times_\z \re$
is the orbit space
$(X\times \re)/\z$ with the action of $\re$ induced by the action
$(x,t)+r=(x, t+r), r\in \re$, on $X\times \re$ and the action of $\z$
considered as a subgroup of $\re$.  We  identify $X$   with the subset 
 $X \times_\z \re$ corresponding to $X\times \{0\}$ in $X \times \re$.
 Note that
 $X\times_\z \re$ is metric compact,
 the action of $\z$ on $X\times_\z \re$ extends the original action
 of $\z$ on $X$,
   $X \times_\z \re= X+[0,1]$ and  for every $r \in \re$ 
 the set $X+r$ is invariant under the action  of $\z$.
We denote by $\pi_X $, $\pi_\re$ and $\pi_\z$  the projections
of $ X \times \re$  to $X$, $\re$ and $X \times_\z \re$ respectively and denote by $\pi$ the map
$\pi : X\times_\z \re \lo S^1=\re/\z$   induced by $\pi_\re$.
The space $X \times_\z \re$, also known  in Topological Dynamics as the mapping torus of $(X,\z)$,  can be  represented as the quotient
of $X \times [0,1]$ where each point $(x,0)$ is identified with $(x+1, 1)$.

\end{section}

\begin{section}
{ Proof of Theorem \ref{borel}}
\label{proof-borel}
 Clearly we may assume $\mdim X < \infty$.
 Fix $\epsilon >0$ and $ d >\mdim X$, and take  an open cover
 $\cal A$ of $X$ such that $\mesh ({\cal A }+z) < \epsilon/2 $
 for every integer $z\in \z$ such that
 $0\leq zd  <n= \ord {\cal  A}$.  Consider Borel's construction $X \times_\z \re$.
 Without loss of generality we may   assume 
  that $\mesh ({\cal A}+r) < \epsilon/2$ for every $r \in \re$ 
  such that
 $0\leq rd  <n+2d$.
Take a finite closed cover $\cal B$ of $X \times [0,1]$
 such that $\ord {\cal B}\leq \ord{\cal A} +1=n+1$
  and 
 $\pi_X(\cal B)$  refines $\cal A$.  Clearly $\cal C=\pi_\z ({\cal B})$ is 
 a closed finite cover of $X \times_\z \re$
 with  $\ord {\cal C} \leq 2\ord{\cal B}=2(n+1)$.  Recall that $X$ is considered as
 a subspace $X \times_\z \re$    and  note
 that every point of $X \times_\z \re$  not belonging to $X$ is covered
 by at most $\ord {\cal B}$ elements of $\cal C$. Also note that we may choose
 $\cal B$ so that 
 that $\mesh \pi({\cal C})$ is as small as we wish. Then we may assume 
 that $\mesh({\cal C}+r) < \epsilon $ for every $r \in \re$ such that $0\leq r d < n+d$.
 
 Recall that $\re$ acts on $X \times_\z \re$, define
 ${\cal D}_i ={\cal C}+(i/n)$ for 
 $0\leq i <  n $   and note that $\mesh( {\cal D}_i  +z) < \epsilon$
 for every $z \in \z$ such that $0\leq z <[n/d]$. Denote $z_i=[n/d]i$.
  Then for the cover
 ${\cal D}=
 {\cal D}_0 \vee ({\cal D}_1 -z_1)\vee \dots \vee ({\cal D}_{n-1}-z_{n-1})$ of
 $X\times_\z \re$ we have that $\mesh({\cal D}+z) < \epsilon$ for every
 $z \in \z$ such that $0\leq z <[n/d]n$.

 Now consider the set $P$ of the points of $X \times_\z \re$ of $\ord> n+1$ with respect to ${\cal C}$.
 Note that $P $ is a closed subset of $X \subset X\times_\z \re$.
 Then $P_i=P+(i/n)$ is the set of the points of $\ord > n+1$ with respect
 to ${\cal D}_i$ and $P_i \subset X+(i/n)$.
 Replace the elements of ${\cal D}_i$  by slightly larger open sets to  get
 an open cover ${\cal D}'_i$ of $X\times_\z \re$ with $\ord {\cal D}'_i =\ord {\cal D}_i$, denote by $P'_i$
 the set of points of $\ord > n+1$ with respect to ${\cal D}'_i$ and
 consider a  canonical map $f_i : X \times_\z \re  \lo K_i=N({\cal D}'_i)$ to a simplicial complex
 $K_i$ (the nerve of ${\cal D}'_i$) for which  the fibers of $f_i$ refine
 ${\cal D}'_i$ and $\dim K_i \leq 2n+1$.  Assuming that  the elements of ${\cal D}'_i$ are sufficiently
 close to the elements of ${\cal D}_i$ and $\mesh \pi({\cal C})$ is small enough
 we can also assume that there are   subcompexes 
 $L_i\subset K_i$ such that $f_i(P'_i+z)\subset  L_i$, 
 $\dim (K_i\setminus L_i) \leq n$,  and $f_i(f_j^{-1}(L_j)+z)\subset K_i\setminus  L_i $
 for every $i$, $j\neq i$ and $z \in \z$.
 
 Consider the map $ f : X\times_\z \re \lo K=K_0 \times \dots \times K_{n-1}$
 defined by $f(x)=(f_i(x+z_i))_{0\leq i < n} $.
   Then the fibers of $ f$ refine
 the cover ${\cal D}'=
 {\cal D}'_0 \vee ({\cal D}'_1-z_1)\vee \dots \vee ({\cal D}'_{n-1}-z_{n-1})$  of $X \times_\z \re$ and
 $f(X \times_\z \re) \subset  M=M_0 \cup  \dots \cup M_{n-1}$
 where $M_i$ is the product of $K_i$ with all  $K_j \setminus L_j$ for  $j\neq i$.
 Note that $\dim M_i \leq 2n+1 +n(n-1)=n^2+n+1$ and hence $\dim M \leq  n^2+n+1$. Thus ${\cal D}'$ can be refined by
 an open cover  of $\ord \leq n^2+n+2$. Recall that 
 $\mesh({\cal D}+z) < \epsilon$ for every
 $z \in \z$ such that $0\leq z <[n/d]n$. Assuming that the elements ${\cal D}'$
 are sufficiently
 close to the elements of ${\cal D}$ we can assume that the latter property
 holds for ${\cal D}'$ as well and we get that $\mdim X \times_\z \re \leq d$ since
 $\frac{n^2+n+2}{[n/d]n}$ goes to $d$ as $n$ goes to $\infty$.
 This implies that $\mdim X \times_\z \re =\mdim X$.

\end{section}

\begin{section}{Kolmogorov's covers}
\label{review-kolmogorov} For a subset $A \subset \re$ and $ r\in \re$  denote $rA =\{ ra : a \in A \}$.
Fix a natural number $n$  and  a real number $\epsilon >0$ and consider 
the collections of closed intervals ${\cal A}_i =\{\frac{\epsilon}{m+1} [{zm+i}, {zm+i +m}]: z \in \z \}$
for $1\leq i \leq m+1$.
Note that $\mesh{ \cal A}_i < \epsilon$ and the collections
${\cal A}_1,  \dots, {\cal A}_{m+1}$ cover $\re$ at least $m$  times. For  a natural $n$ denote
${\cal  A}_i^n =\{ A_1 \times \dots \times A_n: A_1, \dots, A_n\in {\cal A}_i \}$.
Then for $ n\leq m$ the collections ${\cal A}^n_1 \dots {\cal A}^n_{m+1}$ cover $\re^n$ at least
$m-n+1$ times. Indeed, for a point $x=(x_1, \dots, x_n)\in \re^n$ each coordinate of $x$
is covered by  at least $m$ collections  from ${\cal A}_1, \dots, {\cal A}_{m+1}$ and therefore 
there are at least  $m-n+1$ collections ${\cal A}_i$ so that  each of them covers all the coordinates of $x$ and
and hence for each such collection ${\cal A}_i$ the point  $x$ is covered by ${\cal A}_i^n$.
Clearly  ${\cal A}^n_1, \dots, {\cal A}^n_{m+1}$ are collections of disjoint cubes  of $\diam < \epsilon$ with respect to the $l^\infty$-norm in $\re^n$.
We will refer  to ${\cal A}={\cal A}^n_1 \cup \dots \cup {\cal A}^n_{m+1}$ as
Kolmogorov's cover of $\re^n$.

Let $\cal U$ be a finite open cover of a compact metric $X$ with $\ord{ \cal U}\leq n+1$
and let $m$ be any integer such that $n \leq m$.
Consider a canonical map $f : X \lo K =N({\cal U})$ to the nerve  of $\cal U$. Since $K$ is
a finite simplicial complex with $\dim K \leq n$ there is a finite-to-one map 
$g : K \lo \re^n$. 
 Then for a sufficiently small $\epsilon>0$ we have that for Kolmogorov's cover
$\cal A$   the components of the elements 
of $(g\circ f)^{-1}({\cal A})$ will refine the elements of $\cal  U$ and hence by splitting
the elements of $(g\circ f)^{-1}({\cal A}^n_i)$ into finitely many disjoint closed sets and
taking into account that $g(K)$ meets only finitely many elements of $\cal A$ we get
 from each $(g\circ f)^{-1}({\cal A}^n_i)$ a finite family ${\cal F}_i$ of disjoint closed sets
 such that ${\cal F}={\cal F}_1 \cup \dots  \cup {\cal F}_{m+1}$ refines $\cal U$
 and covers $X$ at least  $m-n+1$ times.
  We will refer to $\cal F$ as a Kolmogorov-Ostrand cover.

\end{section}

\begin{section}{Lindenstrauss' level functions}
\label{review-lindenstrauss}

Let $(X, \z)$ be a non-trivial minimal dynamical system on a compact metric $X$,
$U$  a non-empty open set $U$  in $X$ 
and $\phi : X \lo [0,1]$ any map such that $\phi(X \setminus U)=1$ and $\phi^{-1}(0)$
has non-empty interior. Define the random walk on $X$  by stopping at $x \in X$ with 
the probability $1-\phi(x)$ and moving from $x$ to $x-1$ with the probability $\phi(x)$.
Then  Lindenstrauss' level function $\xi : X \lo \re$  is defined at $x \in X$
as the expectation of the number of steps in the random work  starting at $x$. In other words,
$\xi (x)=\phi(x)(1-\phi(x-1))+2\phi(x)\phi(x-1)(1-\phi(x-2))+\dots$
Since the action of $\z$ is minimal, each random walk will eventually hit $\phi^{-1}(0)$ and stop there,
moreover, the number of steps in each random walk is bounded by a number depending only
on the set  $\phi^{-1}(0)$. Thus $\xi$ is well-defined and continuous.
 It is easy to see that $\xi (x+1)=\phi(x+1)(\xi (x)+1)$ for every  $x \in X$ and therefore 
 given a natural number $n$ we have that
 $\xi (x+z)=\xi (x)+z$ for   every integer $-n\leq z \leq  n$ and $x \in X \setminus ((U-n) \cup \dots \cup (U+n))$.
We will refer to $\xi(x)$ as a Lindenstrauss level function determined by $ U$.

\end{section}

\begin{section}{Proof of Theorem \ref{theorem}}
\label{proof-theorem} Let us first present  some auxiliary notations and properties used in the proof.

Let  $(Y, \re)$ be a dynamical system, $A$ a subset of $Y$,  $\cal A$ a collection of subsets of $Y$ and $\alpha, \beta \in \re$ positive numbers.
The subset $A$ is said to be {\bf $(\alpha, \beta)$-small if $\diam( A +r) < \alpha$} for every $r \in [-\beta, \beta]\subset \re$.
The collection  $\cal A$ is said to be {\bf $(\alpha, \beta)$-fine} if  $\mesh({\cal A}+r)< \alpha$ for every $r \in [0, \beta]\subset  \re$.
The collection $\cal A$ is said to be {\bf   $(\alpha, \beta)$-refined} at  a subset $W \subset Y$ if the following two conditions hold:
(condition 1) no element  of ${\cal A}+r$ meets
the closures of both  $W+r_1$ and $W+r_2$ for every $r,r_1, r_2 \in [-\beta, \beta]\subset \re$  with $|r_1-r_2|\geq 1$ and 
(condition 2) if for an element $A$ of ${\cal A}$ the set $A+[-\beta,\beta]$ meets the closure of $W+[-\beta, \beta]$ 
then $\diam (A+r) < \alpha$ for every $r \in [-\beta, \beta]\subset \re$.

\begin{proposition}
\label{prop-1}
Let $(Y, \re)$ be a free dynamical system on a compact metric $Y$, $\omega$ a point in $Y$ and let $\alpha$ and  $\beta $ be positive real numbers.
Then there is an open neighborhood  $W$ of $w$  and an open cover $\cal V$ of $Y$ such that $\ord {\cal V} \leq 3$ and $\cal V$ is $(\alpha, \beta)$-refined at $W$. 

 \end{proposition}
{\bf Proof.}
Note that $L=w+[-3\beta,3\beta]$ is an interval in $Y$ and consider a map
 $g=(g_1,g_2): Y \lo \re^2$ such that $g_1$  embeds $L$ into $\re$ and $g_2(y)$ is the distance to $L$ from $y \in  Y$.
 Then $g^{-1}(a)$  is a singleton for every   $a \in g(L)$ and therefore 
 there is  a fine open cover $\cal O$ of $\re^2$ with $\ord {\cal O}=3$ and 
 a small  neighborhood $W$ of $w$ in $Y$ such that 
 the cover ${\cal V}= g^{-1}(\cal O)$ of $Y$  is  $(\alpha,\beta)$-refined at $W$. Clearly $\ord {\cal V} \leq 3$.
$\black$

\begin{proposition}
\label{prop-2} 
Let $q>2$ be an integer. Then there is  a finite collection ${\cal E}$ of  disjoint closed intervals in $ [0, q) \subset \re$
such that $\cal E$  splits into the union ${\cal  E}={\cal E}_1 \cup \dots \cup {\cal E}_q$  of $q$ disjoint subcollections 
having the property that for every $ t \in \re$  the set  $t +\z\subset \re$  meets at least  $q-2$   subcollections ${\cal E}_i$
(a set meets a collection  if there is a point of the set that is covered by the collection). Moreover,  we may assume
that $\mesh {\cal E}$ is as small as we wish. 
\end{proposition}
{\bf Proof.}
Let $E_1, \dots, E_{q-1}\subset [0,q)$ be $q-1$ disjoint closed intervals of length $>1$. 
Clearly for each $t \in \re$ the set $t +\z$ meets all the intervals $E_i$. Consider a  finite  set 
$A \subset [0,q)$ such that   $1\in A$  and the numbers of $A$ are linearly independent over $\q$
 (the rationals). Define $\phi: \re \lo  \re$ by  
$$\phi(t)=\inf \{ | t+z_1-a_1|  +|t+z_2 -a_2|  : a_1, a_2 \in A, 
z_1,z_2 \in \z, (z_1,a_1)\neq (z_2,a_2) \}.$$ 
Note that $\phi$ is continuous and periodic with period  $1$ and $\phi(t)>0$ for every $t$. Then 
$\sigma=\inf \{ \phi(t) : t \in \re \} >0$. Set  $\Omega$  to  be 
 the open $\sigma/3$-neighborhood  of $A$ and note that  $\Omega$  contains at most one point
of $t+\z$ for every $t \in \re$. Also note  that for every $1\leq i \leq q-1$ the set $E_i \setminus \Omega$ splits into 
a finite collection of disjoint
closed  intervals  and   denote this collection by ${\cal E}_i$. Thus we have that for every $t \in \re$ the set
$t+\z$ meets at least $q-2$ collections from ${\cal E}_1, \dots, {\cal E}_{q-1}$.  Assuming that
$A$ is sufficiently dense in $[0,q)$ we can get $\mesh{\cal E}_i$ as small as we wish for all $i$.
Setting  ${\cal E}_q$ to be empty or any finite collection 
of disjoint  small closed intervals in $[0,q)$ that don't meet
${\cal E}_1, \dots, {\cal  E}_{q-1}$  we get  the collection ${\cal E}={\cal E}_1 \cup\dots\cup {\cal E}_q$
with    the required properties.
$\black$
\\\\
{\bf Proof of Theorem \ref{theorem}.}
Let   $f=(f_1, \dots, f_k):X \lo [0, 1]^k$ be any  map and let $\epsilon>0$ and $\delta>0$
be such that under each $f_i$ the image of every subset of $X$ of $\diam <3\epsilon$ 
is  of $\diam < \delta$ in $[0,1]$. Our  goal is to approximate $f$ by a $\delta$-close map  $\psi$ 
such that  the  fibers of $\psi^\z$ contain at most
$\gamma=\gamma(d)$ points with pairwise distances larger than  $3\epsilon$ where  $\gamma(t)=[k/(k-t]k/(k-t), t<k$.

Since $[\gamma]=[\gamma(t)]$ for every $t>d$ sufficiently close to $d$ 
we can replace $d$ by a slightly larger real number and assume  that $\mdim X < d$.
 Take natural numbers $n$ and $q$ such that
$\mdim X < n/q < d$  and set $m=qk$. Then $n<qd<qk=m$ and $m/(m-n) \leq k/(k-d)$. 
By Theorem \ref{borel} we have $\mdim X \times_\z \re =\mdim X $. Then, assuming that $n$ is large enough, there is 
an open cover ${\cal U}$ of $X \times_\z \re$  such that $\ord {\cal U} \leq n-2$ and $\cal U$ is $(\epsilon, q)$-fine.
Clearly we may assume that $q>2$.

Since the theorem obviously  holds if $X$ is a singleton, we may assume that $(X,\z)$ is  non-trivial. Fix a point $w \in X $ and
 let $l >q$ be a  (sufficiently large) positive  integer which  will be defined later and will depend only on $q$.
By Proposition \ref{prop-1} there is an open cover $\cal V$ of $X\times_\z \re$ and a neighborhood $W$ of $w$ in $X \times_\z \re$ such
 that $\cal V$ is $(\epsilon, 2l)$-refined at $W$. 
 
 Now replacing 
 $\cal U$ by an open cover of $\ord\leq n$  refining ${\cal U}\vee {\cal V}$ 
 we can assume that ${\ord {\cal U}} \leq n$, $\cal U$ is $(\epsilon, q)$-fine
and  $\cal U$ is $(\epsilon,2l)$-refined at $W$. Clearly we can replace $W$ by any smaller neighborhood of $w$ and assume
that $W$ is $(\epsilon, 3l)$-small and  the elements of ${\cal D}_W$ are disjoint where
${\cal D}_W$ is the collection of the closures of $W+z$ for the integers 
 $z\in [-2l,2l]$.
 
Refine $\cal U$ by a  Kolmogorov-Ostrand cover
  ${\cal F}$ of $X \times_\z \re$   such that  
 ${\cal F}$ covers $X \times_\z \re$ at least $m-n$ times and
 ${\cal F}$ splits into ${\cal F}={\cal F}_1 \cup \dots \cup {\cal F}_{m}$ the union of finite  families of
 disjoint closed sets ${\cal F}_i$. Note that $\cal F$ is $(\epsilon,q)$-fine and $(\epsilon, 2l)$-refined at $W$.

 Let $\xi : X \lo \re$ be a Lindenstrauss level function determined by $W $ restricted to $X$.
Denote $W^+ =W +\z \cap [-l,l]$ and $X^-=X \setminus W^+$.
Recall  that 
$\xi(x+z)=\xi(x) +z$ for every  $x \in X^-$ and  an integer $-l\leq z\leq l$.

We need an additional auxiliary notation.
Let  $\cal A$ be  a collection of subsets of $X \times_\z \re$, $\cal B$ a collection of intervals in $\re$. 
For $B\in \cal B$ and $z \in \z$ consider the collection ${\cal A}+B$ restricted to $\xi^{-1}(B+qz)$ and
denote by ${\cal A}\oplus_\xi B$ the union of such collections for all $z\in \z$.   Now  denote
by ${\cal A}\oplus_\xi {\cal B}$ the union
of the collections ${\cal A}\oplus_\xi  B$ for all $ B \in {\cal B}$. Note that ${\cal A}\oplus_\xi {\cal B}$ is
a collection of subsets of $X$.

Consider  a finite  collection ${\cal E}$ of  disjoint closed intervals in $ [0, q) \subset \re$ satisfying  the conclusions of Proposition \ref{prop-2}.
 For $1\leq i \leq k$  define the collection  ${\cal D}_i$  of subsets of $X$ as the union of the collections
${\cal F}_{i}\oplus_\xi {\cal E}_1$, ${\cal F}_{i+k}\oplus_\xi {\cal E}_2$, \dots, ${\cal F}_{i+(q-1)k}\oplus_\xi {\cal E}_q$. 
Note that assuming that $\mesh \cal E$ is small enough we may also assume that ${\cal F}_i^+={\cal F}_i+[-\mesh{\cal E}, \mesh{\cal E}]$
is a collection of disjoint sets and 
 the collection ${\cal F}^+={\cal F}+[-\mesh{\cal E}, \mesh{\cal E}]$ is 
$(\epsilon,q)$-fine and $(\epsilon, 2l)$-refined at $W$  and, as a result,  we get that
 ${\cal D}_i $ is a collection of disjoint  closed sets of $X$ of $\diam < \epsilon$ and
each  element of ${\cal D}_i$ meets at most one element  of ${\cal D}_W$.
 
Define a map $\psi=(\psi_1, \dots, \psi_k) : X\lo [0,1]^k$  so that  for each $i$ the map  $\psi_i$ is $\delta$-close to $f_i$,
$\psi_i$ sends the elements of ${\cal D}_W$ restricted to $X$  and the elements of  ${\cal D}_i$  to singletons in $[0,1]$ and 
$\psi_i$ separates  the elements of ${\cal D}_W$ restricted to $X$  together with  the elements of ${ \cal D}_i$ not  meeting ${\cal D}_W$. 
We will show that the fibers of $\psi^\z$ contain at most  $\gamma$ points with pairwise distances  larger than $  3\epsilon$. 

Denote by $S$ the set of all the pairs   of integers  $(i,j)$  with $0 \leq i \leq l-1$ and $1\leq j \leq k$.
We say that a point $x \in X$ is marked by a pair $(i,j) \in S$ if $x+i$ is covered by ${\cal D}_j$ 
and denote by $S_x \subset S$  the set of  the pairs by which $x$ is marked.
Let us compute the size of $S_x$ for a point $x \in X^-$. Consider a non-negative integer  $z$  such that
$\xi(x) \leq zq< (z+1)q< \xi(x) +l$.  Recall  that $x+(zq-\xi(x))$ is covered by at least $m-n$ collections 
from the family ${\cal F}_1, \dots, {\cal F}_m$ and $\xi(x)+\z$ meets   at least 
$q-2$ collections from ${\cal E}_1, \dots, {\cal E}_q$. 
Then ($\dd$) the point      $x $ is marked by at least $m-n-2k$ pairs $(i,j)$ of $S$ with $zq\leq \xi(x)+i <(z+1)q$.

Indeed, for every ${\cal E}_p$ that meets  $\xi(x) +\z$ pick up $z_p \in \z$ such that $\xi(x)+z_p$ is covered by
${\cal E}_p$.  Denote  $i_p=z_p +zq$ and note that  
$zq \leq \xi(x) +i_p <(z+1)q$  and,  hence,  $0\leq i_p <l$.
 Also note that  different $p$ define different $i_p$ and
for every $1\leq j\leq k$ such that ${\cal F}_{j+(p-1)k}$ covers the point  $x+(zq-\xi(x))$ we have 
that the collection ${\cal D}_j$ covers $x+i_p$, and therefore $x$ is marked by the pair $(i_p,j)$ in $S$.
Thus if $\xi(x)+\z$ meets all the collections ${\cal E}_1, \dots, {\cal E}_q$  
the number of   pairs $(i,j) \in S$ with $zq\leq \xi(x)+i <(z+1)q$ marking $x$ will be 
at least the number of times $x +(zq-\xi(x))$ is covered 
by the collections ${\cal F}_1, \dots, {\cal F}_m$,  which is
at least $m-n$.
Each time  $\xi(x)+\z$ misses a collection from ${\cal E}_1, \dots, {\cal E}_q$  reduces the above estimate by at most $k$. 
Since $\xi(x)+\z$ can miss at most two collections from ${\cal E}_1, \dots, {\cal E}_q$ we 
arrive at the required estimate $m-n-2k$ in ($\dd$).

Thus, by ($\dd$), the point $x$ is marked by at least $(\frac{l}{q}-3)(m-n-2k)$ pairs of $S$. Since  $m-n\geq m(k-d)/k=q(k-d)$ we have 
$|S_x|\geq (\frac {l}{q}-3)(m-n-2k) \geq ((l/q)-3)(q(k-d)-2k) =l(k-d)\Delta$
where $\Delta =(1-\frac{3q}{l})(1-\frac{2k}{(k-d)q})$.
Note that we may take $q$ sufficiently large and  then $l$ sufficiently large with respect to $q$ and assume
that $\Delta < 1$ and $\Delta$ is as close to $1$ as we wish.

Now we will show that there is no set contained in a fiber of $\psi^\z$ and  containing more than  $\gamma$ points with pairwise
distances larger than $3\epsilon$.  Aiming at  a contradiction assume that such  a  set  $\Gamma \subset  X$  exists with $|\Gamma |=[\gamma]+1$.
First  note  that    $\Gamma$  can contain at most  one point  in $X\setminus X^-$ since each $\psi_i$ separates ${\cal D}_W$ restricted to $X$
and $\mesh({\cal D}_W)\leq \epsilon$.

Now assume that $\Gamma$  contains a point $x \in X\setminus X^-$ and
a point $y \in X^-$. Take   any pair $(i,j) \in S_y$
and an element $D \in {\cal D}_j$ that contains $y+i$.  Since $\psi^\z_j$ does not separate $x$ and $y$ we get that  $D$ meets 
the element of ${\cal D}_W$ containing $x+i$, and ($\dd\dd$) the latter is impossible  because ${\cal F}^+$ is
$(\epsilon, 2l)$-refined at $W$ and $W$ is $(\epsilon,3l)$-small that  implies that 
 $x$ and $y$ are $2\epsilon$-close. 
 
 Indeed, $D$ is contained in an element of ${\cal F}^+  +t$ for some $0\leq t < q$. Then $y$ is covered by an element
 of ${\cal F}^+ + t-i$ that intersects the element of ${\cal D}_W$ containing  $x$, and the facts that 
$(\epsilon, 2l)$-refined at $W$ and $W$ is $(\epsilon,3l)$-small  yield the conclusion of ($\dd\dd$).

Thus, by ($\dd\dd$), we get  that $\Gamma \subset X^-$. Now we can assume that $\Delta > \gamma/|\Gamma|$ and  get 
that  $\sum_{x \in \Gamma} |S_x | >|\Gamma| l(k-d)(\gamma/|\Gamma|)=lk [k/(k-d)]$. Then, since $|S|=lk$,
there is a pair $(i_*,j_*)$ in $S$  and $\Gamma_*\subset \Gamma$ such that 
$|\Gamma_*|= [k/(k-d)]+1$ and every point in $\Gamma_*$ is marked by $(i_*,j_*)$.
Thus $\Gamma_* +i_*$ is covered by ${\cal D}_{j_*}$.
Then  ($\dd\dd\dd$) there is an element $D$ of ${\cal D}_{j_*}$ containing $\Gamma_*+i_*$.

 Indeed, if $\Gamma_*+i_*$
meets two elements of ${\cal D}_{j_*}$ then those elements cannot touch the same element
of ${\cal D}_W$ (because the points of $\Gamma_*$ are $3\epsilon$-distant and  
${\cal F}^+$ is $(\epsilon, 2l)$-refined  at $W$) and, hence will be separated by $\psi_{j_*}^\z$.  
Thus  ($\dd\dd\dd$) holds.

Then, by ($\dd\dd\dd$), we have that
 $\xi(\Gamma_*) $ is contained in an element (interval) of ${\cal E} +qz$ for some non-negative integer $z$.  
 Since $\mesh{\cal E}$ is small there is an integer $0\leq q_*\leq q-1$ such that $\xi(\Gamma_* -q_* ) \subset[ qz,qz+2]$.
 Denote by $S_*$ the set of  pairs $(i,j)$ such that $0\leq i\leq q-3, 1\leq j \leq k$,
 and   say that a point $x\in \Gamma_* -q_*$ is marked by $(i,j) \in S_*$ if
 $x+i$ is covered by ${\cal D}_j$.
 Since ${\cal F}$ covers $X$ at least $m-n$ times we get by an  argument similar to the one
applied above that each point of $\Gamma_*-q_*$ is marked by at least $m-n-6k\geq q(k-d)-6k=q(k-d)\Delta_*$ pairs in $S_*$
where $\Delta_*=1-\frac{6k}{q(k-d)}$. Note that  $(k/(k-d))/|\Gamma_*| < 1$ and  assume that $q$ is such  that $(k/(k-d))/|\Gamma_*|<\Delta_*$.
Then, since
$|S_*|=k(q-2)$,  there are two distinct points $x$ and $y$  of $\Gamma_*-q_*$ marked by the same pair $(i,j)$ of $S_*$. 
Again by the above reasoning  $x+i$ and $y+i$ cannot lie 
in different elements of ${\cal D}_j$. 
Thus $x+i$ and $y+i$  lie in the same element $D$ of ${\cal D}_j$.  Then, since
${\cal F}^+$ is $(\epsilon, q)$-fine, 
this  implies that   the points $ x+q_*$ and $y+q_*$ of  $\Gamma_*$  are $\epsilon$-close. 

Indeed, $D$ is contained in an element of ${\cal F}^+ +t+i$ with $t\in [0,2)$.    Then  $x$ and $y$ are
contained in an element of  ${\cal F}^+ +t$, and   $x+q_*$ and $y+q_*$ are contained in an element
of ${\cal F}^+ +t+q_*$. Hence $x+q_*$ and $y+q_*$ are $\epsilon$-close since $0\leq t+q_* < q$  and
${\cal F}^+$ is $(\epsilon, q)$-fine.

Thus  the assumption $|\Gamma| > \gamma$ leads to a contradiction and hence
  $\psi$ is the desired approximation of  $f$. Then the theorem follows by a standard Baire category  argument.
$\black$

\end{section}

 Department of Mathematics\\
Ben Gurion University of the Negev\\
P.O.B. 653\\
Be'er Sheva 84105, ISRAEL  \\
 mlevine@math.bgu.ac.il\\\\ 
  
\end{document}